\documentclass[graybox]{svmult}

\usepackage{mathptmx}       

\usepackage{helvet}         
\usepackage{courier}        

\usepackage{graphicx}              

\usepackage{siunitx}
\usepackage{amsbsy}
\usepackage{amsmath,bm}
\usepackage{amssymb}
\usepackage{amsfonts}

\usepackage{multirow}

\usepackage[bottom]{footmisc}
\usepackage{cite}

\usepackage{url}

\begin{document}

\title*{Solution of time-harmonic Maxwell's equations
by a domain decomposition method based on PML
transmission conditions}
\titlerunning{PML based transmission conditions for Maxwell's equations in DDM}
\author{Sahar Borzooei 
\and
Victorita Dolean 
\and
Pierre-Henri Tournier 
\and
Claire Migliaccio 
}
\institute{Sahar Borzooei, Victorita Dolean \at Côte d'Azur University, CNRS, LJAD,France,
\email{Sahar.Borzooei@univ-cotedazur.fr, work@victoritadolean.com}
  \and Pierre-Henri Tournier \at Sorbonne University, CNRS, LJLL, France,
  \email{tournier@ljll.math.upmc.fr}
  \and Claire Migliaccio \at Côte d'Azur University, CNRS, LEAT, France, \email{Claire.Migliaccio@univ-cotedazur.fr}}

%
%
\maketitle

\abstract*{ Numerical discretization of the large-scale Maxwell's equations leads to an ill-conditioned linear system that is challenging to solve. The key requirement for successive solutions of this linear system is to choose an efficient solver. In this work we use Perfectly Matched Layers (PML) to increase this efficiency.
PML have been widely
used to truncate numerical simulations of wave
equations due to improving the accuracy of the solution instead of using absorbing boundary conditions (ABCs).
Here, we will develop an efficient solver by providing an alternative use of PML as transmission conditions at the interfaces between subdomains in our domain decomposition method. 
We perform a series of numerical simulations on our model and will assess the convergence rate and accuracy of our solutions compared to the situation where absorbing boundary conditions are chosen as transmission conditions.}

\abstract{Numerical discretization of the large-scale Maxwell's equations leads to an ill-conditioned linear system that is challenging to solve. The key requirement for successive solutions of this linear system is to choose an efficient solver. In this work we use Perfectly Matched Layers (PML) to increase this efficiency.
PML have been widely
used to truncate numerical simulations of wave
equations due to improving the accuracy of the solution instead of using absorbing boundary conditions (ABCs).
Here, we will develop an efficient solver by providing an alternative use of PML as transmission conditions at the interfaces between subdomains in our domain decomposition method. 
We solve Maxwell's equations and assess the convergence rate of our solutions compared to the situation where absorbing boundary conditions are chosen as transmission conditions.}

\section{Introduction}
\label{Borzooei:Intro}
Maxwell's equations need to be solved in many applications, such as medical imaging or electromagnetic compatibility. The Finite Element
Method (FEM) is widely used for numerical modeling of these problems due to its ability to handle complex geometrical configurations. Finite element discretization of these frequency-domain wave problems leads
to an ill-conditioned linear system with large number of unknowns. To solve this system, the efficiency of direct solvers is limited at larger scales due to scalability problems in memory and computing time. Besides, Krylov subspace iterative solvers have shown slow convergence. An alternative method that tackles the convergence problem of iterative solvers is the domain decomposition method (DDM). The method relies on a division of the computational domain into smaller
subdomains, leading to subproblems of smaller sizes manageable by direct solvers. One perfect candidate introduced in \cite{Borzooei-b1} and then improved in \cite{Borzooei-b15} is the domain decomposition preconditioner which proved to be very robust in large scale computations.
 However, designing an efficient domain decomposition preconditioner is still challenging for such a system. \\
In this paper, we present an efficient PML-based Schwarz-type domain decomposition preconditioner with overlapping subdomains. The convergence rate of Schwarz methods highly depends on the transmission condition on the interfaces between subdomains. Thus, carefully designed transmission conditions play a critical role in the efficiency of the solver. To decrease undesired numerical reflections, one usually adds a PML layer along the boundaries that extend a definite area to the infinity, representing an unbounded volume,  and absorbs almost all incident waves, regardless of angle of incidence, so that the waves decay exponentially in magnitude into the PML medium \cite{Borzooei-1,Borzooei-4}. Specifically, using PML is essential for simulating unbounded systems such as infinitely long waveguides or an isolated structure in an infinite vacuum region.  While the use of PMLs as boundary conditions when solving a problem in open space is quite common, less things are known about their use as transmission conditions within a domain decomposition algorithm. We propose to assess the performance of a one-level domain decomposition algorithm where the transmission conditions at the boundaries between subdomains are PML conditions, providing a better approximation to the transparent boundary operator. We will investigate the convergence properties and compare them with the more common impedance transmission conditions.
In a previous work, PML have been used successfully as transmission conditions in domain decomposition methods in geophysical applications modeled by the Helmholtz equation in \cite{Borzooei-tournier}. The paper is organized as follows. In Section 2, we present mathematical model including Maxwell's equations, PML formulation with different stretching functions as well as its implementation. Then, the DDM with PML-based transmission operators is introduced.
In Section 3, some numerical examples are presented to analyze the performance of the proposed domain decomposition algorithm. Finally, conclusion is written in Section 4.

\section{Mathematical Model}
\label{Borzooei:sec_2}
Let us consider the computational domain $\Omega \subset \mathbb{R}^3$ to be a homogeneous dielectric medium of complex-valued electric permittivity $ \varepsilon_{\sigma} $ and electrical conductivity $ \sigma >0 $. Let  $\mu_0 $ be the permeability of free space and \textbf{n} be the unit outward normal to the boundaries $\partial\Omega$. $\omega$ is the angular frequency  and $c $ is the wave speed. In the frequency domain, the electric field $\mathbf{\xi}(\mathbf{x},\mathbf{t})= \Re(\mathbf{E(x)}e^{i \omega \mathbf{t}})$ has harmonic dependence on time of angular frequency $\omega$, where $\mathbf{E(x)}$ is its complex amplitude depending on the space variable $\mathbf{x}$.  Hence $\mathbf{E(x)}$ is a solution to the following second order time-harmonic Maxwell’s equation
\begin{equation}
\label{Borzooei:eq_max}
\begin{array}{rcl}
   \nabla\times(\nabla\times \mathbf{E})- \omega^2 \varepsilon_{\sigma}\mu_0\mathbf{E}~ = ~f~~~~~~~~~~~~\text{in }~ \Omega.~~~~~~~~~~~~
\end{array}
\end{equation}
Let us denote the boundary of the global domain by $\partial\Omega$ where  Robin condition $(\nabla\times \mathbf{E})\times \textbf{n}+ i\frac{\omega}{c} \textbf{n} \times (\mathbf{E}\times \textbf{n})~= 0$ is imposed~\cite{Borzooei-sah}. The Robin or impedance boundary condition (Imp BCs) is a standard first order approximation to the far field Sommerfeld radiation condition enabling the description of the wave behavior in a bounded domain, while the physical domain is not bounded.
The finite element discretization of equation~\eqref{Borzooei:eq_max} is written as the following linear system

\begin{equation}
\label{Borzooei:eq_lin}
A\mathbf{u}=\mathbf{b}.
\end{equation}
 
\subsection{PML Formulation}
\label{borzooei:subsec_1}
To solve a partial differential equation (PDE) numerically, the computational domain has to be truncated without introducing reflections. The first attempt in this regard is absorbing boundary conditions (ABCs). The first order ABC as regular choice is Robin condition that was mentioned earlier. Due to the limited accuracy of this method, PML were introduced by Berenger\cite{Borzooei-1} as a better alternative.  
PML provide non-reflecting boundaries so that the numerical solution converges exponentially to the exact solution in the computational domain as the thickness of the layer increases. \\
PML implementation is done by stretching cartesian coordinates so that stretching is defined in a layer surrounding the computational domain \cite{Borzooei-1} and Dirichlet boundary condition can be imposed at the end of the PML layer. 
In this regard, we assume the boundaries of the computational domain to be aligned with the coordinate axes. For simplicity we will focus on truncating the problem in the $x$ direction. Let us suppose that the PML layer extends from the boundary of our domain of interest $x=a$ until $x = a^{*}$. The coordinate mapping in $x$ direction is:
\begin{equation}
\label{Borzooei:eq_pml}
\frac{\partial}{\partial x_{pml}} \rightarrow \frac{1}{1 - \frac{i}{\omega}\sigma(x)} \frac{\partial }{ \partial x}~~~~, 
\sigma(x) = \begin{cases}
\begin{array}{ll}
~~~0 ~ ~~~~~~~~~~ \text{if } x<a  ~~\\
>0 ~~~~~~ \text{if } a < x < a^{*}.~~
\end{array}
\end{cases}
\end{equation}

In the PML region where $\sigma(x)>0$, The oscillating solutions turn into exponentially
decaying ones. In the rest of the region where $\sigma(x)=0$, the wave equation is unchanged and
the solution is unchanged. In this paper we have studied two different stretching functions $\sigma(x)$ as following
\begin{equation}
\label{Borzooei:eq_loss}
\begin{array}{ll}
\sigma_{-1}(x) = \frac{1}{a^{*}-x} ~~~~~~~~(a)~~~~~ , ~~
   \sigma_{-2}(x)= \frac{2}{(a^{*}-x)^{2}}~~~~~~~(b)
\end{array}
\end{equation}

To truncate our computational region with a PML layer in other directions, we just need to do the same transformations to get $\frac{\partial}{\partial y_{pml}}$ and $\frac{\partial}{\partial z_{pml}}$. At the corners of the computational
cell, we will have PML regions along two or three directions simultaneously, but it will not generate any problem. \\
Implementing this mapping in the three dimensional domain requires a slight further generalization of equation~\eqref{Borzooei:eq_max}, resulting in the following definition of the curl operator to be used in the variational formulation:

\begin{equation}
\label{Borzooei:eq_curl}
 \nabla_{pml}\times \mathbf{E} =
        \left[ {\begin{array}{c}
      \frac{\partial{\mathbf{E}_z}}{\partial{y_{pml}}}-\frac{\partial{\mathbf{E}_y}}{\partial{z_{pml}}} \\
    
      \frac{\partial{\mathbf{E}_x}}{\partial{z_{pml}}}-  \frac{\partial{\mathbf{E}_z}}{\partial{x_{pml}}} \\
     
    \frac{\partial{\mathbf{E}_y}}{\partial{x_{pml}}}- \frac{\partial{\mathbf{E}_x}}{\partial{y_{pml}}} \\
    \end{array} } \right]
\end{equation}
\subsection{Domain Decomposition Preconditioner}
 To solve our large and ill conditioned linear system \eqref{Borzooei:eq_lin}
, the use of a robust and efficient preconditioner is necessary in a Krylov iterative solver (GMRES)\cite{Borzooei-4}. A preconditioner $M^{-1}$ is a linear operator that approximates the inverse of matrix $\textbf{A}$
whose cost of the associated matrix-vector product is much cheaper than solving the original linear system. In this regard, we employ right preconditioning to solve \eqref{Borzooei:eq_lin} that will give us:
\begin{equation}
\label{Borzooei:right}
AM^{-1}\mathbf{y}=\mathbf{f}, ~~~~\textit{where} ~~~ \mathbf{u}=M^{-1}\mathbf{y}
\end{equation}
This right preconditioned system benefits from a residual that is preconditioner independent compared to the left-preconditioned variant.

 As an overlapping Schwarz method, the optimized restricted additive Schwarz (ORAS) domain decomposition preconditioner is chosen here
\begin{equation}
M_{ORAS}^{-1}= \sum_{s=1}^{N_{sub}} R_s^T D_s A_s^{-1} R_s  \label{Borooei:oras}
\end{equation}
where $N_{sub}$ is the number of overlapping subdomains $\Omega_s$ into which the domain $\Omega$ is decomposed.  Here, matrices $A_s$ stem from the discretisation of local boundary value problems on $\Omega_s$ with transmission conditions at the subdomain interfaces. Let $N$ be an ordered set of the unknowns of the whole domain and let $N = \bigcup_{s=1}^{N_{sub}}N_s$ be its decomposition into the (nondisjoint) ordered subsets corresponding to the different (overlapping) subdomains $\Omega_s$. Matrix $R_s$ is the restriction matrix from $\Omega$ to subdomain $\Omega_s$; it is a $ N_s \times N$ Boolean matrix. $R_s^T$ is then the extension matrix from subdomain $\Omega_s$ to $\Omega $. $D_s$ is a $N_s \times N_s $ diagonal matrix that gives a discrete partition of unity, ie, $ \sum_{s=1}^{N_{sub}} R_s^T D_s R_s = I$. \\

The convergence rate of this method highly depends on the choice of transmission conditions between the subdomains\cite{Borzooei-6}. The optimal convergence is obtained by imposing the Dirichlet-to Neumann (DtN) map related to the complementary of each subdomain \cite{Borzooei-v2}, \cite{Borzooei-v2b}. However, since the cost of computing the exact DtN is prohibitive,
low-order absorbing boundary conditions (ABCs) to approximate the DtN have been developed. Nonetheless, these methods have limited accuracy, which led to developing domain decomposition strategies with high order transmission conditions \cite{Borzooei-v3}. But the problem with high order transmission conditions is the difficulty of their implementation. A good approximation of ABCs in terms of providing better convergence rate and easy implementation would be to use PML on the interface boundaries of the cuboid-shaped subdomains \cite{Borzooei-v1}, \cite{Borzooei-v1c}, that is what we consider here. In this purpose a PML layer is added in each direction in the overlap region. Note that the width of the overlap has to be larger than the PML layer for a good transmission of the data between subdomains.

\section{Numerical Results}
\begin{figure}
\centering
\scalebox{0.35}{\includegraphics{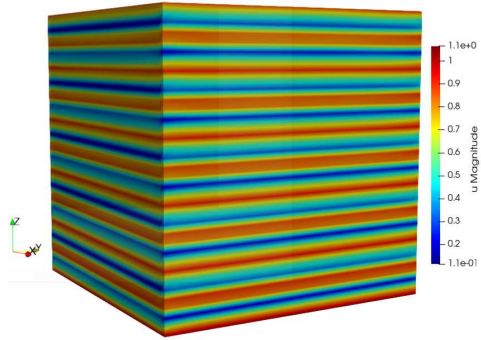}}
\caption{Plane wave propagation.}
\label{Borzooei:fig_1}       
\end{figure}

The performance of the proposed PML-based preconditioner for Maxwell's equations is studied in a 3D homogeneous domain $\Omega$, while length of the domain in each direction is $10~m $. We have excited the $z = 0$ surface with plane wave incident term $e^{(-ikz)}$, where $k = \frac{2\pi}{\lambda}$, with propagation in $+z$ direction shown in the Fig. \ref{Borzooei:fig_1}.
 The convergence rate is studied while there is PML or Impedance as global boundary conditions (BCs) or interface conditions (ICs), which leads to four different situations reported in table~\ref{Borzooei:tab_1}. The finite element discretization is done for the first order edge elements for two different frequencies. Let $\#$DoF represents the number of degrees of freedom. For $f = 0.5~Hz$, we have $\#$DoF =511775 and for $f = 1~Hz$, we have $\#$DoF = 2098100. The global domain is decomposed into N = 100 number of cuboid-shaped subdomains, that PMl layer is set along their interfaces with the length $L_{pmli}$. \\ 
  In table~\ref{Borzooei:tab_1}, simulations are done for the $\sigma_{-2}$ stretching function and PML length on the interfaces is shown with $L_{pmli} = 8h$ where $h = \frac{\lambda}{n_{\lambda}}$ is the mesh size and $n_{\lambda}$ is the number of points per wavelength. In all simulations, PML length on the global boundary is $L_{pml}=2\lambda$ and the overlapping layers between subdomains is changed from 2 to 8 layers in four steps. Looking at the table~\ref{Borzooei:tab_1}, for $f=0.5$ and considering 8 number of overlapping layers to be larger than $L_{pmli}$, we can see while we have PML BCs and Imp ICs, number of iterations is 21, but with PML BCs and PML ICs, this number decrease to 16. It is while with Imp on BCs and ICs, number of iterations is 26. In this table, $\bullet$ means that, solution has not converged for 200 number of iterations.

\begin{table}
\centering
\caption{ function $\sigma_{-2}$. $n_{\lambda}$ = 5,$L_{pml} = 2\lambda$, $~L_{pmli} = 8h$ , $c =1$, $N = 100$ is number of subdomains }
\label{Borzooei:tab_1}  
		\begin{tabular}{cc|cccc|cccc}
		 &&&f =0.5 &&&~~~~~~~~~~f =1&&\\
		\hline
		BCs & ICs&2&4&6&8&2&4&6&8\\
		\hline
		\hline
		Imp & Imp  &29&24& 23&26&34&27& 25&25\\
		Imp & PML  &75&30& 22&28&$\bullet$&50&29 &23\\
		PML & Imp &23&21&20 &21&28&23& 21&20 \\
		PML & PML &65&25& 18&16&$\bullet$&43&25 &19\\
		\hline
	\end{tabular}
\end{table}

To see the influence of the $L_{pmli}$, we did the simulation with smaller PML layer on the subdomains, mentioned in Table~\ref{Borzooei:tab_2}, only for the case with PML as BCs and ICs. Comparing with the equivalent row in the Table~\ref{Borzooei:tab_1} for 6 and 8 overlaping layers, we see that with larger length of PML on subdomains we have better convergence. Although the rate of convergence has become better for lower overlapping layers, with smaller PML length, due to the better data transmission between subdomains. Comparing the results for the use of stretching function $\sigma_{-1}$  instead of $\sigma_{-2}$ is mentioned in the Table~\ref{Borzooei:tab_3}. For 2 number of overlapping layers, better number of iterations is seen with $\sigma_{-1}$, however for higher number of overlapping layers, $\sigma_{-2}$ results in better convergence. 
\begin{table}[h]
\centering
\caption{ function $\sigma_{-2}$. $n_{\lambda}$= 5, $L_{pml} = 2\lambda$, $L_{pmli}= 4h$ , $c =1$, $N = 100$ }
\label{Borzooei:tab_2}  
		\begin{tabular}{cc|cccc|cccc}
		 &&&f =0.5 &&&~~~~~~~~~~f =1&&\\
		\hline
		 BCs & ICs&2&4&6&8&2&4&6&8\\
		\hline
		\hline
		PML & PML &59&22& 18&17&183&35&30 &29\\
		\hline
	\end{tabular}
\end{table}	
	
\begin{table}
\centering
\caption{function $\sigma_{-1}$. $n_{\lambda}$ = 5, $L_{pml} = 2\lambda$, $L_{pmli} = 8h$ , $c =1$, $N = 100$ }
\label{Borzooei:tab_3}  
		\begin{tabular}{cc|cccc|cccc}
		 &&&f =0.5 &&&~~~~~~~~~~f =1&&\\
		\hline
		BCs & ICs&2&4&6&8&2&4&6&8\\
		\hline
		\hline
		PML & Imp &24&20& 18&21&29&24&22 &21 \\
		PML & PML &46&28&19 &17&66&35&26&21\\
		\hline
	\end{tabular}
\end{table}

The performance of the proposed preconditioner in a heterogeneous domain is studied in Tables~\ref{Borzooei:tab_4} and \ref{Borzooei:tab_5}. Here, we have defined a medium with two values of $\varepsilon_r$ with the dimension of $6.3~m$ in X and Y directions and $2.5~m$ in Z direction inside the free space computational domain. In this experiment, rhs is chosen as a random value, $L_{pmli} = 4h$, $f=1~Hz$, stretching function is chosen as $\sigma_{-2}$ in Table~\ref{Borzooei:tab_4} and $\sigma_{-1}$ in Table~\ref{Borzooei:tab_5}. Results show, increasing value of $\varepsilon_r$ increase number of iterations, but with PML interface conditions we can have faster convergence. Comparing two tables, better performance is obtained by $\sigma_{-2}$ stretching function.  Here we have considered maximum number of iterations as 600. In the results, - means problem is not solved due to memory limitation.

\begin{table}
\centering
\caption{ function $\sigma_{-2}$. $n_{\lambda}$ = 5, $L_{pml} = 2\lambda$, $~f=1.0~Hz$ , $c =1$, $N = 100$.}
\label{Borzooei:tab_4}  
		\begin{tabular}{cc|cccc|cccc}
		 &&&~~$\varepsilon_r =4$ &&&&~~~$\varepsilon_r =5$\\
		\hline
		BCs & ICs&2&4&6&8&2&4&6&8\\
		\hline
		\hline
		PML & Imp&262 & 207&181&168&586& 425&375&326\\
		PML & PML& $\bullet$&249 &199&160&$\bullet$& 514&440&311\\
		\hline
	\end{tabular}
\end{table}

\begin{table}
\centering
\caption{ function $\sigma_{-1}$. $n_{\lambda}$ = 5, $L_{pml} = 2\lambda$, $~f=1.0~Hz$ , $c =1$, $N = 100$.}
\label{Borzooei:tab_5}  
		\begin{tabular}{cc|cccc|cccc}
		 &&&~~~$\varepsilon_r =4$ &&&&~~~~$\varepsilon_r =5$\\
		\hline
		BCs & ICs&2&4&6&8&2&4&6&8\\
		\hline
		\hline
		PML & Imp&264&207&183&-&587&427 &374&-\\
		PML & PML&411&221&201&-&$\bullet$& 402&353&-\\
		\hline
	\end{tabular}
\end{table}

Results are obtained on the Université Côte d’Azur's High-Performance Computing (HPC) center. In this HPC center,
cluster is composed of 48 CPU computing nodes, including
32 nodes with Dual Intel Xeon Gold processor, providing 40
cores per node and 192 GB of memory and 16 nodes with 2
AMD Epyc processors, providing 32 cores per node and 256
GB of memory. 

\section{Conclusions} In this work, we have developed a numerical model for an accurate and fast simulation of Maxwell's equations. To achieve this goal, the PML layer is implemented as physical boundaries and as transmission conditions in domain decomposition preconditioner for a three dimensional domain. A better convergence rate is achieved with PML layer, compared to Impedance interface conditions. Numerical results shows that the performance of the PML depends on a well chosen stretching function and length of the PML. This work is a preliminary study that was inspired by  a similar work done for Helmholtz equations\cite{Borzooei-sahDD} where the results were very encouraging.
More investigations can be done in next works, like evaluating performance of the PML as interface conditions for higher order edge elements or in a heterogeneous domain. Note that PML have some limitations, for the time being it has been applied only along the straight interfaces but variants for circular boundaries exist that can be further explored in the context of other applications. 

\begin{acknowledgement}
This project has received funding from the European Union’s Horizon 2020 research and innovation programme under grant agreement No 847581 and is co-funded by the Région Provence-Alpes-Côte d'Azur and IDEX $UCA^{JEDI}$. This work was supported by the French government, through the $UCA^{JEDI}$ Investments in the Future project managed by the National Research Agency (ANR) under reference number ANR-15-IDEX-01. The authors are grateful to the OPAL infrastructure from Université Côte d’Azur and the Université Côte d’Azur’s Center for High-Performance Computing for providing resources and support.
\end{acknowledgement}



\end{document}